\documentclass[unnumsec,webpdf,contemporary,large]{oup-authoring-template}%
\usepackage{caption}

\graphicspath{{Fig/}}

\theoremstyle{thmstyleone}%
\newtheorem{theorem}{Theorem}%  meant for continuous numbers
\theoremstyle{thmstyletwo}%
\theoremstyle{thmstylethree}%

\begin{document}

\journaltitle{PNAS Nexus}
\DOI{DOI HERE}
\copyrightyear{2022}
\pubyear{2022}
\access{Advance Access Publication Date: Day Month Year}
\appnotes{Paper}

\firstpage{1}

%\subtitle{Subject Section}

\title[Soft Matter Complexity and Computing]{From Navier-Stokes millennium-prize problem to soft matter computing}

\author[a,b,$\ast$]{Saksham Sharma}
\author[c]{Giulia Marcucci}

\authormark{Sharma et al.}

\address[a]{\orgdiv{Pembroke College}, \orgname{University of Cambridge}, \orgaddress{\street{Pembroke Street}, \postcode{CB2 1RF}, \state{Cambridge}, \country{UK}}}
\address[b]{\orgdiv{Department of Chemical Engineering}, \orgname{University of Cambridge}, \orgaddress{\street{Phillipa Fawcet Dr.}, \postcode{CB3 0AS}, \state{Cambridge}, \country{UK}}}
\address[c]{\orgdiv{Apoha Ltd., London}, \orgname{}, \orgaddress{\street{Acklam Rd.}, \postcode{W10 5JJ}, \state{London}, \country{UK}}}

\corresp[$\ast$]{To whom correspondence should be addressed: \href{email:ss2531@cam.ac.uk}{ss2531@cam.ac.uk}}

\received{Date}{0}{Year}
\accepted{Date}{0}{Year}
%\editor{Associate Editor: Name}

\abstract{Clay Mathematical Institute, in the year 2000, formulated a list of seven unsolved mathematical problems which might influence and direct the course of the research in the 21st century. Navier-Stokes regularity problem is one of the seven problems and has not been solved till date. These equations govern the motion of the viscous fluids, such as liquids, gases, gels, polymers. The present article is a spotlight on the NSE problem: a review of the past and current efforts to solve it, and our changing perspectives and understanding of these equations. The major debate is to either prove that the solutions to Navier-Stokes equations (NSE) are smooth or to prove that the solutions reach singularity at a finite-time (blowup); given some initial conditions. While in the past, the focus was on finding bounds for the smooth solutions, recently, there has been growing interest in finding solutions that blowup. The current article places the Turing-completeness of Euler equations (a subclass of NSE) into the centre-stage, by discussing the recent developments that are pointing towards that the idea that the fluid equations should be viewed as a ``computer program''. Given this is the case, the article then argues to see ``continuum fluid" from a complexity notion, and thus, as a Turing machine. This notion connects the discipline of fluid mechanics to disciplines like, \textit{neuromorphic computing, reservoir computing} and potentially leading to emerging new discipline, \textit{soft matter complexity, computing, and learning}, each of which has been reviewed in brief. The understanding that is unfolding from this article is that the unsolved classical problems in fluid dynamics, such as turbulence, ``continuum hypothesis'', ``NSE regularity'', need a change in approach from reductionist mathematical style (such as, PDE analysis, chaos theory) to a theoretical computational one (Turing machine, computational complexity, machine learning).}

%\boxedtext{
% Research reports require a significance statement of between 50 and 120 words. Abbreviations are permitted, but citations cannot be included. If required, un-comment this element in the template to include. The heading is included automatically.
%}
\maketitle

Dating back to the times of Aristotle, water has been considered as one of the fundamental elements of nature (along with earth, air, and fire). This view, though, was demolished during the Chemical Revolution of the late eighteenth century, when it was accepted that water is, instead, a compound made of the elements hydrogen and oxygen \citep{iswaterh20}. Ever since then, our understanding of water has progressed enormously, from knowing the constituent elements to their subatomic (electrons, protons, neutrons). Not only water, our scientific understanding of `fluids' in general (liquid, gas, gels, polymers) - materials that flow when deformed - has improved by leaps and bounds; we now know to a greater microscopic details (inter-atomic forces, packing density) - experimentally and theoretically - the dynamical behaviour of such fluids \citep{dave_soft_rev}. \newline

\begin{figure*}
  \centering
  \includegraphics[width=1.0\textwidth]{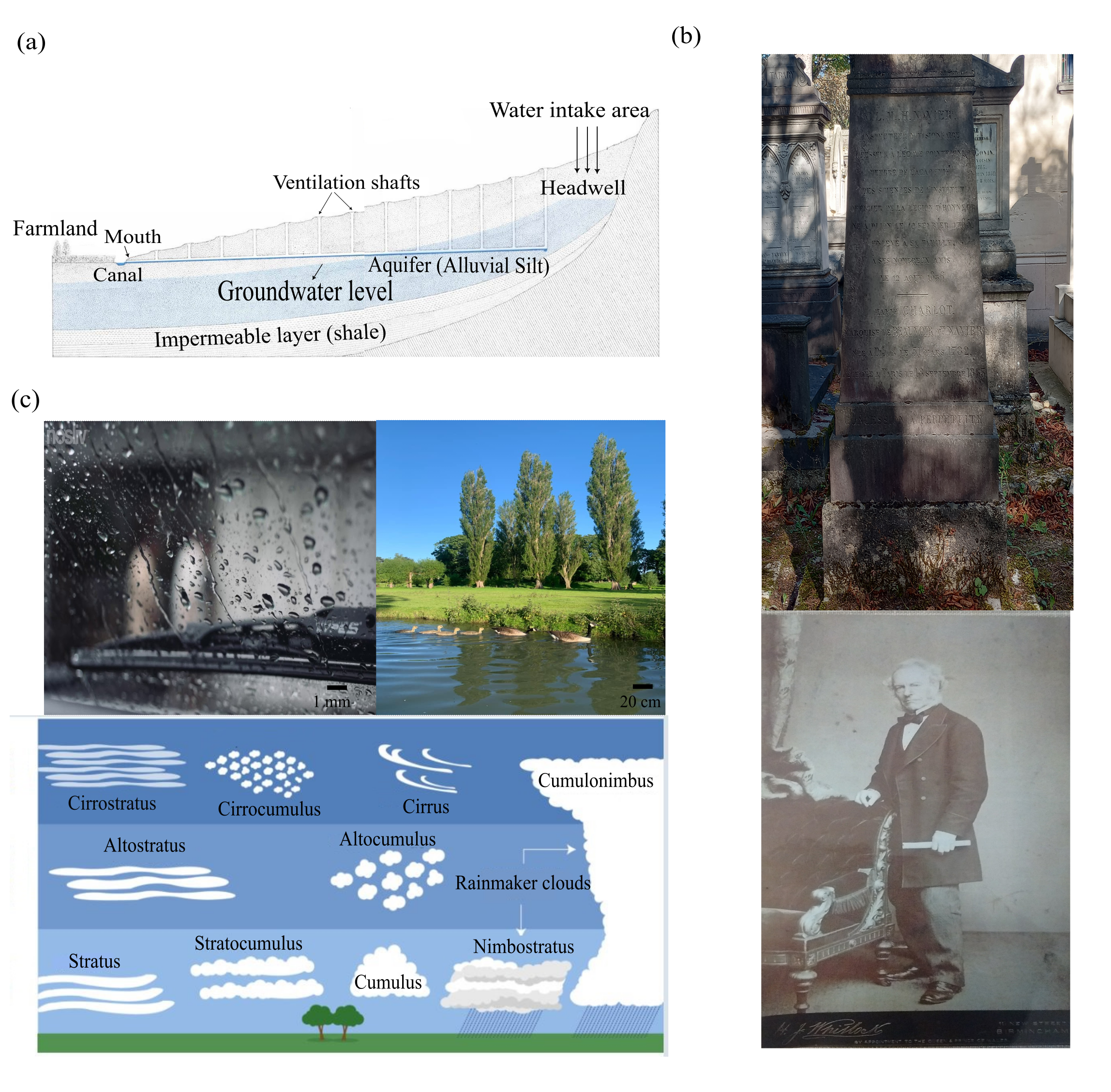}

 \caption{A collage on the Navier-Stokes equations. (a) Large-scale constructions, called `quanats' in the ancient Iran built by Persian people around 1st millennium BCE to transport water from a headwell to the farmland through the vertical shafts. Adapted from \citep{wulff1968qanats} ; (b) French engineer and physicist Claud-Louis Navier's graveyard (top) at Père-Lachaise, division 50 in Paris (France) photographed by S.S. An undated photograph of Irish physicist Sir George Stokes (bottom) (photograph taken during \href{https://stokes200.weebly.com/}{\textcolor{blue}{Stokes 200 centenary conference}} at the Pembroke College); (c) Photographs of the water drops on the window pane of a car (top-left), photo of a mother duck being trailed by four baby ducks at the Grantchester Meadows in Cambridge (top-right), schematic of a variety of cloud shapes (bottom, adapted from: \href{https://scied.ucar.edu/learning-zone/clouds/cloud-types}{\textcolor{blue}{here}}).}
  \label{fig:fig1}
\end{figure*}

Understanding of fluids has always inspired engineering advances: hydraulics, in the form of force pumps developed by Greek scientists; `quanats', constructed in ancient Iran by the Persians (Fig. \ref{fig:fig1}(a)); principle of buoyancy for floating objects, discovered by Archimedes - have been known for millennia. However, the scientific revolution in the seventeenth century, beginning with the works of Galileo, Newton, and Leibniz, led to the discovery of physical laws that explain the motion of accelerated objects. It was Leonard Euler, who set aside experiments, and applied the celebrated Newton's second law of motion to fluid parcels/packets (also called ``continua'') - an emergent collection of those becomes a fluid \citep{euler1757principes}. In the words of Clifford Truesdell, ``looking within the interior moving fluid, where neither eyes nor experiment may reach, he (Euler) called upon the imagination, fancy, and invention which Swift could find neither in music nor in mathematic'' \citep{truesdell1960program}. In a continental-wide effort, the analysis of fluids by Euler (a Swiss mathematician working in Russia) was, subsequently, improved by a French physicist Claude-Louis Navier (Fig. \ref{fig:fig1}(b)), by incorporating nonlinear forces within the fluid that arises because of the stretching of fluid parcels, and finally an Irish physicist George Gabriel Stokes (Fig. \ref{fig:fig1}(b)) combined all the efforts together to derive equations for a viscous fluid, which we now call the ``Navier-Stokes equations'' (NSE) \citep{stokes2007theories}. \newline

Ever since their formulation, NSE (with some added force components) have given this world a mathematical formalism to understand tides, winds, ocean flows, sea level rise, motion of vehicles, and almost all engineering marvels that encounter fluid flow. Next time when you behold pearls of water drops running down the window pane of your car \citep{podgorski2001corners}, ducklings following their mothers while swimming in a river \citep{yuan2021wave}, or wave-like clouds formed (fluctus) with a clear blue sky in the background (Fig. \ref{fig:fig1}(c)), recall that the NSE can explain these phenomena to a reasonable extent \citep{miyazaki2001method}. However, a complete understanding of NSE is lacking, in the sense that its solutions are not always trustworthy. For certain cases, they might be physically reasonable (smooth), for other cases, they might be irregular and yield singularity (blowup). There is no consensus on either of the two possibilities, which is precisely the core of millennium-prized problem formulated by Clay Mathematical Institute in the beginning of the 21st century. \newline

This article starts with summarising the statement of the Navier-Stokes regularity problem in Section \ref{NSR}. Then the article discusses the efforts for past one century by mathematicians in finding analytical solutions and bounds for smooth solutions to the equations, in Section \ref{sec:efforts}. Section \ref{new_scheme} discusses a series of papers for the past one decade that are pointing towards the fundamentally `computational' nature of the equations. Given this is the case, the article argues a position in Section \ref{consequences} that \textit{fluids} should be viewed as \textit{computational entities}, and hence subject to the computational treatment experimentally. This idea leads one to treating fluids using ideas from the discipline of neuromorphic computing (section \ref{neuromorphic}). The article then concludes in the section \ref{outlook} with the possible consequences of this idea on the NSE regularity problem and soft matter research in the next few decades. \newline

\section{Navier-Stokes regularity problem}\label{NSR}
Navier-Stokes equations, are Newton's laws of motion to model the motion of a fluid, in two or three spatial dimensions ($\mathcal{R}^2 \textrm{or} \, \mathcal{R}^3$), and single time $t$ dimension. The velocity vector field $\vec{u}(x,t) = (u_{i}(x,t))_{1 \leq i \leq n} \in \mathcal{R}^n$ and scalar pressure field $p(x,t) \in \mathcal{R}$ are the unknowns in the equations. Given an assumption of incompressibility of the fluid, these equations are given by 
\begin{subequations} \label{eq:ns}
\begin{gather}
    \partial_t \vec{u}(t,x) = \nu \Delta \vec{u} - \sum_{i=1}^3 u_i \partial_i \vec{u} -  \vec{\nabla} p + \vec{f} \\
    \nabla.\vec{u}=\sum^3_{i=1} \partial_i u_i = 0 
\end{gather}
\end{subequations}
where $\nu$ is the kinematic viscosity (units: $m^2 s^{-1}$) and $\vec{f}$ is the external body force. The initial conditions of the velocity vector field are given as 
\begin{equation} \label{eq:init}
    \vec{u}(0,x)=\vec{u}_0 (x)
\end{equation}
where \eqref{eq:ns} and \eqref{eq:init} are defined for $t \geq 0$ and $x \in \mathcal{R}^n$. If $\nu=0$, then the fluid theoretically has no viscosity (inviscid fluid), and the equations are termed \textit{Euler equations}. The solutions to either of these equations are considered physically reasonable if 
\begin{subequations} \label{eq:phys}
\begin{gather}
    p,u \in C^\infty (\mathcal{R}^n \times [0,\infty)) \\
    \int_{\mathcal{R}^n} |u(x,t)|^2 dx < C 
\end{gather}
\end{subequations}
for all $t \geq 0$. Here, $C^\infty$ refers to the functions $p,u$ being smooth and infinitely differentiable, and the kinematic energy of the fluid is bounded above by a constant $C$. \newline

\noindent The original problem statement postulated at the dawn of 21st century by Clay Mathematical Institute \citep{fefferman} is to prove either of the following statements:
\begin{enumerate}
    \item For $\nu>0$ and $n=3$, prove that for the smooth, divergence-free, physically reasonable initial vector field $u_0(x)$, and $f(x,t)=0$, there exists smooth functions $p(x,t),u_{i}(x,t)$ on $\mathcal{R}^3 \times [0,\infty]$ that satisfy \eqref{eq:ns} - \eqref{eq:phys}. 
    \item For $\nu>0$ and $n=3$, prove that for the smooth, divergence-free, physically reasonable initial vector field $u_0(x)$, and $f(x,t)=0$, there exists no solutions $p(x,t),u_{i}(x,t)$ on $\mathcal{R}^3 \times [0,\infty]$ that satisfy \eqref{eq:ns} - \eqref{eq:phys}. 
\end{enumerate}
It is to be noted that this problem statement for $\nu=0$ is still unsolved, however, Euler equations are not listed as the prize problems. Furthermore, while the \eqref{eq:ns}-\eqref{eq:phys} are defined for spatial domains extended to infinity, the equations can also be reformulated for a spatially periodic case. Such a case is also enlisted in the clay prize problem. \newline

\section{Efforts to solve the problem}
\label{sec:efforts}
Attempts to find solutions to the Navier-Stokes equations (NSE) goes back to early 20th century, long before the problem, stated above, was formally formulated. In 1911, \cite{oseen1911formules} resolved a paradox posed by George G. Stokes known as Stokes' paradox. Stokes' paradox stated that there are \textit{no} bounded (finite) solutions to the Stokes equations (Reynolds number is negligibly small) for the case of a infinitely long cylinder or a disk in 2D falling through a fluid. Oseen resolved the paradox by taking into account the inertial term ($u_{i}\partial_{i}\vec{u}$) and a singular point force (called \textit{stokeslet}). Stokeslet allows the velocity to grow as the logarithm of the radial distance, and thus incorporated as the inner boundary condition to yield bounded solutions $\vec{u}(x,t)$ outside the cylinder.  In 1934, Leray \citep{leray1934mouvement} found a set of solutions (called \textit{weak solutions}) $u_{0}$ belonging to $L^{2}$ norm\footnote{An $L^2$ norm or the Euclidean norm of a function gives the length of a vector $\mathbf{x}=(x_{1},x_{2},...,x_{n})$ as: $||\mathbf{x}||_{2}=\sqrt{x_{1}^2+x_{2}^2+..+x_{n}^2}$}, for a manifold $\mathcal{R}^n$ such that the weak solution in addition to energy inequality \textit{shall be unique} if a strong solution exists (or the weak solution is regular enough to be a strong solution). However, this dream has not been achieved till date because weak solutions have not been proved yet to be unique, thus regular, in nature. The physical meaning of such a \textit{partial regularity} is that the solutions $\vec{u}(x,t)$ match analytically well with the smooth solutions only in the large regions of spacetime, which excludes the regions of null sets (measurable sets of measure zero) \footnote{Cantor sets are one kind of null sets.} corresponding to time $t$ and tiny length scale $x$. Such null sets are finite sets of arbitrary small length, and are usually disconnected, or in other words, discrete. Hence, it is possible to have Leray weak solutions which would be different at different times $t$ (only because of unknown regularity in the null sets), leading one to have non-unique general class of solutions. \newline

In 1993, Constantin,  Fefferman \cite{constantin1993direction} proved that, for given constants $\Omega$ and $\rho$, if the angle between unit vorticity vectors at positions $x$ and $y$ is $\varphi(x,y,t)$, then the magnitude of the vorticity $\omega$ at both locations is above $\Omega$ if
\begin{equation}
    \left| \sin \varphi(x,y,t) \right| \leq \frac{|x-y|}{\rho}  
\end{equation}
such that the vortex tubes are formed in an antiparallel fashion osculating (observed physically and numerically in experiments) each other. The vorticity vector is regular only in $L^{2}$ norm and not until $L^{\infty}$ norm. As a result, this well-behaved solution of vorticity vector (in the sense of parallel or antiparallel alignment) exists only because of the viscous mechanism, leading to \textit{partial regularity} within the nonlinear stretching mechanism. \newline
\noindent One type of scaling of NSE is given as 
\begin{equation}
    v(x,t) \mapsto v_{\lambda}(x,t)=\lambda v (\lambda x, \lambda^2 t)
\end{equation}
and there are various partial regularity results known for this. One of the most prominent one is called Ladyzhenskaya-Prodi-Serrin conditions (\citep{ladyzhenskaya1957existence}, \citep{prodi1959teorema}, \citep{serrin1961interior}), which states that if the Leray weak solution lies in $L^p_t L^q_x$ with $2/p + 3/q \leq 1$, then the solution is unique and smooth in positive time. The endpoint in the above equality for $p=\infty$ and $q=3$, given for $L^{\infty}_t L^3_x$, was proved by  Escauriaza-Seregin-\u Sver\'ak \cite{escauriaza2003img}. It should be noted that while the natural scaling of weak solutions has a regularity for the condition, $2/p + 3/q = 3/2$, the energy equality holds in NSE with an additional regularity for the condition, $2/p + 3/q = 5/4$ \citep{shinbrot1974energy}. This means that there is a gap open between the natural scaling of the equations and the kinetic energy, which can lead to non-uniqueness of weak solutions. In fact, this nonuniqueness was hinted by Jia-\u Sver\'ak \citep{jia2014local}, where the authors proved that in the class $L_t^\infty L_x^{3,\infty}$, Leray solutions are nonunique if a certain spectral assumption holds for a linearised Navier-Stokes operator. In 2019, Buckmaster and Vicol \citep{buckmaster2019nonuniqueness} proved a stronger result that the weak solutions to the 3D NSE, $v \in C_{t}^{0} H_{x}^{\beta}$ (where $H^{\beta}$ is the $L^{2}$ space with the regularity index $\beta$), are nonunique. The nonuniqueness refers to ill-posedness and existence of infinitely many solutions so that a clear statement about the nature of the solution can not be made. The idea used in the proof of this work is the convex integration method \citep{de2015h} and the authors expect that these tools might in the future establish nonuniqueness of Leray weak solutions. \newline

\begin{figure*}

  \centering
  \includegraphics[width=1.0\textwidth]{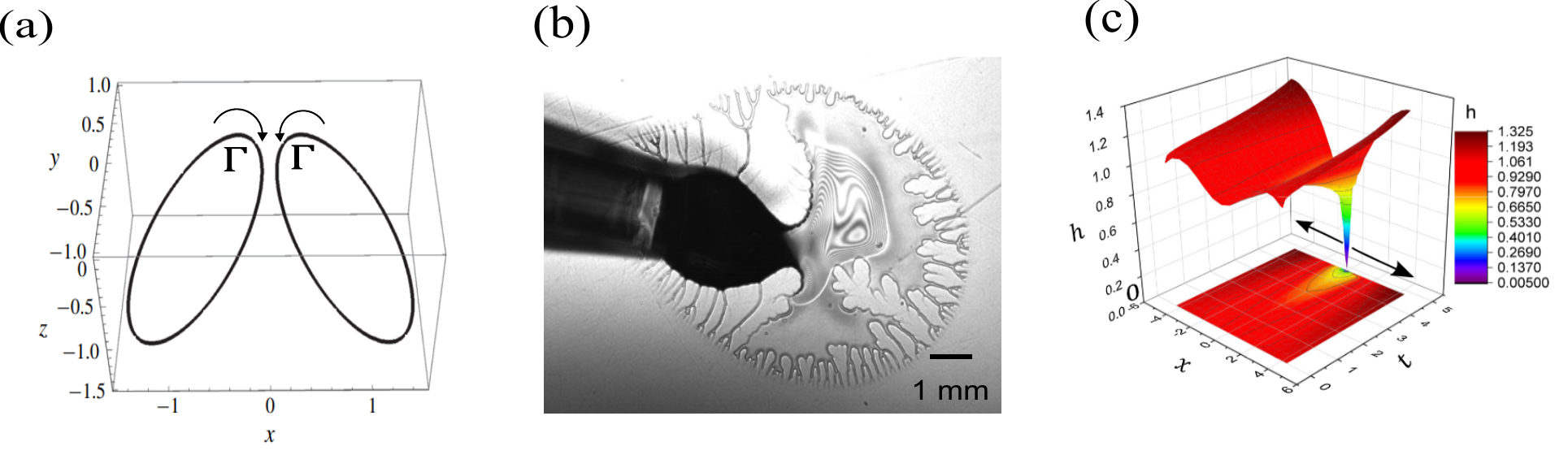}
 \caption{Examples from the fluid dynamics literature where the physical systems attain finite-time singularity. (a) Pyramid vortex reconnection observed by Kimura and Moffatt \citep{moffatt2019towards}, (b) A receding thin-film (on a polyethylene substrate) of a viscoelastic sticky fluid secreted by carnivorous pitcher plants \textit{N. rafflesiana} forms slender filaments which are of fractal-like shapes at the periphery; taken from \citep{kang2021sticky}, (c) A liquid thin-film of height lower than a critical value (which depends on surface tension and interaction forces between liquid and the substrate) becomes unstable and yields singularity ($h=0$) at a finite-time \citep{sharma2022}.}
  \label{fig:fig2}
\end{figure*}

While the results shown above are obtained by functional analysts for past one century, the NSE question has also percolated the community of applied mathematicians. Kimura and Moffatt \citep{moffatt2019towards} showed that two vortex rings, of finite radii and counter-rotating circulation $\Gamma$ (Fig. \ref{fig:fig2}(a)), when approaching each other can exhibit a finite-time singularity, which is physical (not mathematical) in nature, when the key length scale of the system becomes smaller than the intermolecular scale, making the \textit{continuum hypothesis} \footnote{A fundamental physical assumption in the derivation of NSE, continuum hypothesis states that the the effect of statistical and thermal fluctuations because of the motion of molecules in a fluid are not \textit{necessary} to be taken into account while deriving continuum laws of motion, given that the length scale of system is at least 10 times greater than the mean free path of molecules in the continua. See \citep[see p. 6]{batchelor2000introduction}} becomes invalid. In fact, such finite-time physical singularities are fairly ubiquitous in household phenomena such as: the pinch-off of pendant drop (Fig. \ref{fig:fig2}(b).i); fingering patterns in a Hele-Shaw flow cell consisting of parallel plates separated by a distance (Fig. \ref{fig:fig2}(b).ii); eddies formed at an edge of a wedge-shaped container driven by the rotation of a cylinder on the top of the container (Fig. \ref{fig:fig2}(b).iii); slender-like filaments formed by a receding drop of biological fluid secreted by pitcher plants (Fig. \ref{fig:fig2}(b).iv); finite-time singularity formed by an unstable thin-film (Fig. \ref{fig:fig2}(b).v). It is to be noted that the case of vortex reconnection has also been considered by functional analysts to check if the NSE solutions are smooth and regular in nature or not. It turns out that the vortex reconnection has two stages. The first stage is when the strong coherent vortices come together and stretch to form intense narrow vortex tubes which approach at various angles by violating the smoothness assumption on the gradient of direction $\nabla \zeta$. The second stage is when the viscosity plays an important role by having an upper bound
\begin{equation}
\nu^2 \int |\omega| |\nabla \zeta|^2 dx dt \leq E_0
\end{equation}
on the vorticity. As a result, the vortex tubes become antiparallel, reconnect, and finally decay. While this picture has been numerically shown in several studies \citep{meiron1989numerical}, there is no rigorous proof on the regularity of solutions in this problem. 

\section{A new scheme to look at the problem}
\label{new_scheme}
In Section \ref{sec:efforts}, numerous efforts to either prove uniqueness or hint towards nonuniqueness of the solutions to NSE were outlined. This section, specifically, highlights a new scheme to solve the problem by using ideas from theoretical computer science. \newline

\begin{figure*}
  \centering
  \includegraphics[width=1.0\textwidth]{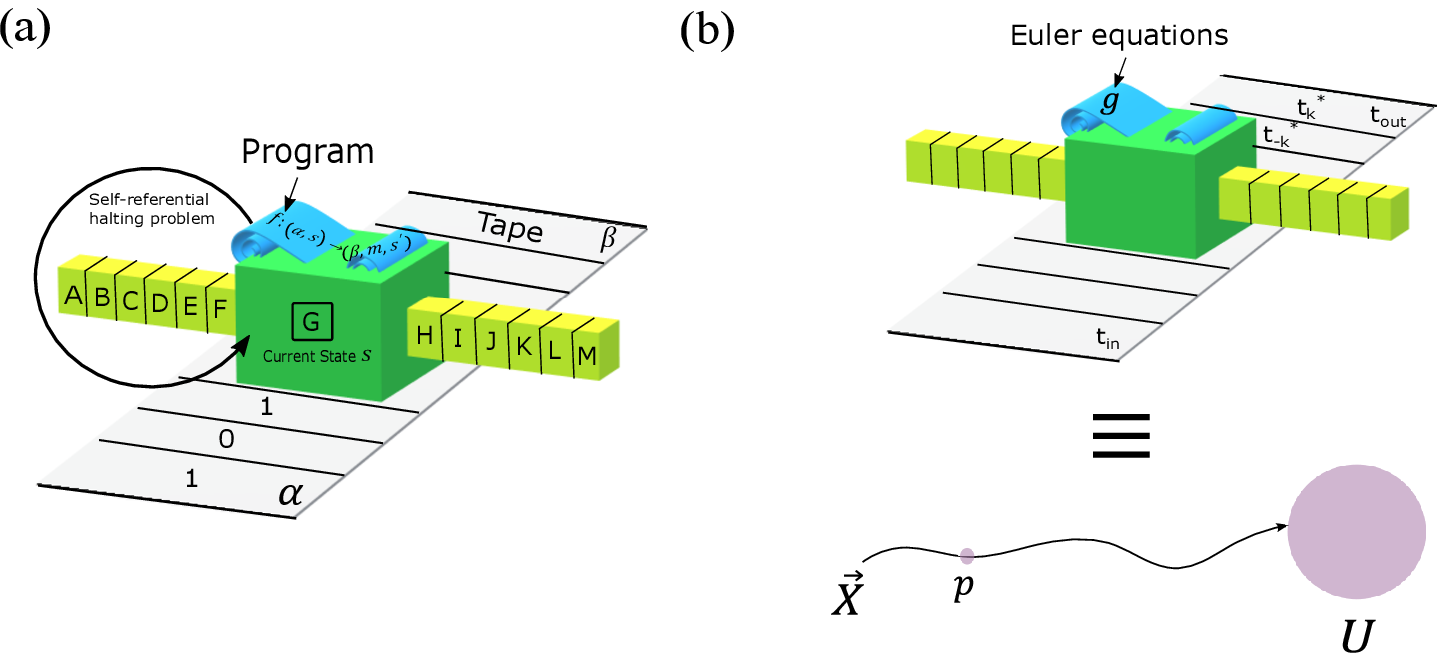}
 \caption{Fluid equations in a computational complexity perspective. (a) Schematic of a Turing machine on a tape with input $\alpha$, program/function $f$, and output $\beta$, (b) Euler equations as embeddable in a Turing machine and suffering from the \textit{Halting problem} is equivalent to the undecidability of a flow field $\vec{X}$ passing through the point $p$ and entering an open set $U$.}
  \label{fig:fig3}
\end{figure*}

\noindent In 2015, Terence Tao \citep{tao2016finite} introduced a model given as
\begin{equation}
\begin{gathered} \label{eq:dyadic_var}
\partial_t X_n = - \lambda^{2n\alpha} X_n + 1_{n-1=n(t)}\lambda^{n-1}X^{2}_{n-1} \\ - 1_{n=n(t)}\lambda^{n=n(t)}\lambda^{n}X_{n}X_{n+1}
\end{gathered}
\end{equation}
where $n: [0,T_{*})\to\mathcal{Z}$ is a piecewise constant function that describes the mode pair ($X_{n(t)},X_{n(t)+1}$) allowed to interact at a given time $t$, $X_{n}$ is the energy of fluid at scale $n-1$ being diffused to scale $n$. A system of ODEs were constructed to \textit{simulate} \eqref{eq:dyadic_var} by using a sequence of ``quadratic circuits" connected in series fashion. Each circuit is composed of ``quadratic logic gates" which simulate a certain form of basic quadratic non-linear interaction. As shown in Fig. \ref{fig:fig3}(c), the gates transfer energy from one mode to another with ``amplifier" and ``rotor" gates. Tao argued that it might be possible that such gates are Turing complete. After 6 years, in 2021, Cardona, Miranda, Peralta-Salas, Presas \citep{cardona2021constructing} proved that it indeed is the case that Euler equations are Turing complete. The theorem that the authors proved states that
\begin{theorem}
There exists a Riemannian manifold $(M,g)$ whose Euler flow is Turing-complete: the halting problem for any given Turing machine is equivalent to the Euler flow for a certain initial condition associated to that machine entering a certain region.
\end{theorem}
where the authors constructed solutions for the steady Euler flow (without viscosity) on a Riemannian 3-sphere $\mathcal{S}^{3}$ that are Turing-complete. Here, Turing completeness means that for given points on the sphere, the problem of bounding the dynamical trajectory of those points, as the equation evolves, is an undecidable problem. Following this work, the authors extended the analysis to a standard Euclidean 3-dimensional space. The tools used to prove such a result are contact topology, symplectic geometry, h-principle, generalised shift maps, and Cantor sets. The key step in the construction of their proof is translating the problem in Fluid Dynamics as a problem in Geometry by using a key master point on establishing a correspondence between stationary Euler flows (Beltrami fields) and Reeb vector fields - a former result by Etnyre, Ghrist and also envisaged formerly by Dennis Sullivan (2022 Abel Prize winner \citep{papadopoulos2022dennis}). Such a mirror is used in the proof to construct a Reeb vector field in dimension three. This construction of the Reeb vector field reminds Poincare's section and promotes a 2D to 3D construction proving that a certain map in the plan can be seen as time-1-map of the flow of a Reeb vector field. The 2D construction builds on former work of Cristopher Moore in 1991 \citep{moore1991generalized} on whether hydrodynamics is able to perform computations, so the construction in this work answers Moore's question. An intuitive sense of how the proof is derived can be given as follows. For a given string of symbols of the alphabet $(t^{*}_{-k},....,t^{*}_{k})$ and a Turing machine $T$ with an input tape $t$, the authors showed that it is possible to construct explicitly a point $p$ and an open set $U$ in $\mathcal{S}^{3}$ such that when the velocity field $X$, which is the solution of inviscid Navier-Stokes equations, passes through $p$, it also intersects $U$ if and only if $T$ halts with the output tape at positions corresponding to symbols $(t^{*}_{-k},....,t^{*}_{k})$ (Fig. \ref{fig:fig3}(d)). This work marks as an important step in towards treating the general solution of different subclasses of NSE as embeddable in a Turing machine. The authors found that if the initial vector field $\vec{X}$ for Euler equations is given as an input in the NSE, the resulting velocity field can only simulate a finite number of steps, and thus, the Turing completeness of Navier-Stokes equations would still remain unproven. \newline

The results presented above are only for manifolds which are of higher dimensions, and thus, are not scalar invariant. Proving finite-time blowup is an issue for such manifolds, as at finer scales, the metric structure becomes flat and the equivalence between Euler flow and arbitrary vector fields does not hold true. The key next step here is to find solutions to Euler equations on $\mathcal{R}^{3}$ which is scalar invariant, as the original Clay problem (\eqref{eq:ns}) states. Tao discussed in MSRI seminar that certain type of hybrid manifolds which exhibit scalar invariance in one region and universality in other region might potentially be a way to move forward to prove Euler flows Turing-complete in $\mathcal{R}^{3}$ be importantly noted that Turing-completeness of Euler equations in $\mathcal{R}^{3}$ is not a straightforward path as it seems, as once it is established that the equations are Turing-complete, the proof of finite-time blowup will require lots of care in the proof. \newline

As discussed in Tao \citep{tao2016finite}, proving the Turing-completeness of the Euler or Navier-Stokes equations is the first step towards building, designing, and quantifying the design of logic gates that can be assembled to simulate the classical laws of conservation of energy and helicity that the fluid equations will obey. After making sure that the physics is right and design of logic gates is rigorous enough, such a ``fluid program'' can be shown to blowup for a specific set of initial data, and the debate of regularity versus blowup can be given more concrete answer, for a given choice of initial data.

\section{Consequences to the fluid dynamics research}
\label{consequences}
As discussed by George Batchelor \citep{batchelor2000introduction}, the mathematical basis for the continuum treatment of liquids in motion is incomplete and thus, it is difficult to deduce the properties of a hypothetical continuous medium that moves the same way as a real fluid. Very recently, there is a growing evidence that `fluids' do not much show the `regular' or symmetric behaviour, at certain regimes where the flow is turbulent or the length scale is close to the molecular scale. We recently argued that there lacks a fundamental framework which connects the physics of fluids at two different scales \cite{sharma2022complexity}. If the analysis of fluids is restricted only to one continuum scale without incorporating the underlying molecular contributions (\cite{bandak2022dissipation}), then the Navier-Stokes equations yield singular solutions in an unpredictable fashion. To solve this problem, the topic of ``emergence" within complexity sciences field comes to our rescue as it celebrates the idea that physics at multiple scales can be correlated or conceptually connected to each other. One framework built off complexity view of fluid mechanics suggested by us (\cite{sharma2022complexity}) is the use of cantor sets (Lindenmayer systems) as a computational framework to derive continuum level properties from molecular level. These systems consists of certain (recursive) rules being applied at one scale to yield the physics at another scale.   \newline

If one accepts the use of L-systems or related computational frameworks in analysing fluids at multiple scales, then another corollary from this result is to start accepting that fluids, in fact, (\citep{wheeler2018information}) consist of \textit{information} (bits) at their core. This changes the way we see fluids; instead of fluids being modelled mathematically, they can be computationally programmed. Such a program can be both experimental (\cite{adamatzky2019brief, brun2022fluid}) and virtual (\cite{sharma2022complexity, cardona2021constructing}) in nature. \newline

Given that the fundamental nature of fluid dynamical behaviour is computational in nature, the question that immediately comes up is: is it digital or analog? To answer this question, it is important to resource two important results. First is that the the fluid dynamical equations are chaotic in nature, as discussed extensively by Doering \& Gibbons \cite[see Ch. 5]{doering1995applied}. Second is that the chaotic behaviour of a physical system can be encoded in \cite{siegelmann1995computation} in analog shift maps which are more powerful than Turing machine, and computes similar to how a neural network and an analog machine would compute. These two results can be combined to argue that a fluid can be programmed to behave like an analog computer, and perform computations similar to a neural network. This result leads one to explore new avenues of research, particularly in \textit{neuromorphic and unconventional computing} communities, as the tools from these fields can be borrowed to program the analog behaviour of fluids, soft matter, or active matter. 
\section{Neuromorphic computing}
\label{neuromorphic}

In the increasing community of scientists working on analog and unconventional computing, neuromorphic computing (NMC) is the branch that has been gaining the most attention since the beginning of the fourth industrial revolution %[https://link.springer.com/chapter/10.1007/978-3-319-05624-1_1].
Its main objective is to overcome the challenge of engineering machines based on brain-like algorithms, namely, capable of solving problems following a holistic approach, a peculiarity exclusive to the human brain. \newline

The way NMC addresses this challenge involves both cutting-edge software architectures and novel hardware paradigms, in a never-ending interchange. Indeed, while it is possible to trace the path that led the scientific community from the need for efficient big-data analysis to the development of artificial neural networks (ANNs) \citep{fernando2003pattern}, nowadays separating the development of software and hardware in NMC is becoming an impossible task. For this reason, modeling fluid dynamics using NSE and analyzing the computing power of a liquid state machine (LSM) can be seen as two sides of the same coin, a concept that justifies the presence of this section in this paper, and on which the authors are actively working and will show original results in forthcoming publications. \newline

From an algorithm perspective, the two most important architectures of NMC are spiking neural networks (SNNs) \citep{fernando2003pattern}, and reservoir computing (RC) \citep{fernando2003pattern}. SNNs are ANNs based on time-dependent synapses. At variance with an ensemble of perceptrons (basic units of ANNs with output defined as $y=f(\sum_i w_i x_i + b)$, where $f$ is the activation function, $x_i$ are the nodes, $w_i$ the corresponding trainable weights, and $b$ is the bias term), the neurons of a SNN are allowed to transmit information only when the activation potential crosses its threshold. \newline

RC represents a subclass of recurrent neural networks (RNNs) with training performed only at the readout stage. In other words, RC can be considered as a RNN where the total ensemble of hidden layers is flattened into a sparse highly-dimensional layer with fixed dynamics (i.e., a black box). It benefits from the RNN's Turing completeness in most of its architectures, excluded the limit ones, like the extreme learning machine (ELM) which is a RC feedforward version.
The most widespread RC paradigms are echo state networks (ESNs), LSM, ELM, and all RC quantum versions, outside the purpose of this paper. In terms of equations, we can generalize their models to
\begin{equation}
\begin{array}{rcl}
\mathbf{y}(t_j) & = & \mathbf{W}^{out}\mathbf{x}(t_j) + \mathbf{b},\\
\mathbf{x}(t_j) & = & f\left(\mathbf{W}^{s}\mathbf{u}(t_j) + \mathbf{W}\mathbf{x}(t_{j-1})\right),
\end{array}
\label{eq:RC}    
\end{equation}
with $\mathbf{y}(t_j)$ the network readout at time $t_j$ satisfying the condition
\begin{equation}
||\mathbf{y}(t_j)-\mathbf{y}_{T}(t_j)||<\epsilon \;\; \forall j=1,...,N,
\label{eq:RC_training}
\end{equation}
where $\mathbf{y}_{T}$ is the target output and $\epsilon\lesssim 0$. In Eqs. (\ref{eq:RC}), $\mathbf{W}^{out}$ is the matrix of trainable weights obtained from Eq. (\ref{eq:RC_training}) as 
\begin{equation}
\mathbf{W}^{out} = \left(\mathbf{y}_{T}(t_1)-\mathbf{b},...,\mathbf{y}_{T}(t_N)-\mathbf{b}\right)\cdot\left(\mathbf{x}(t_1),...,\mathbf{x}(t_N)\right)^{\dagger},
\label{eq:Wout}
\end{equation}
where $\dagger$ stands for the Moore-Penrose pseudo-inverse matrix; $\mathbf{x}(t_j)$ and $\mathbf{u}(t_j)$ are the vectors of nodes in the output and input layers at time $t_j$, respectively; $\mathbf{b}$ is the bias vector; $f$ is the activation function; $\mathbf{W}^{s}$ and $\mathbf{W}$ are the matrices of the non-trainable weights of the reservoir.\newline

The ELM is obtained from Eqs. (\ref{eq:RC}) when $\mathbf{W}=\mathbf{0}$, the limit case characterized by absence of recurrence. It has already been shown that the nonlinear Schr\"odinger equation can model an ELM \citep{marcucci2020theory}, a fact that suggests that also a reservoir governed by the NSE is able to perform interpolation, classification, and logic gates. Indeed, the nonlinear Schr\"odinger equation can be derived from the NSE in deep-water approximation \citep{peregrine1983water}.
More precisely, in Marcucci \textit{et al.} (2020) \citep{marcucci2020theory}, the authors demonstrated how a continuous-wave laser beam propagating in a self-focusing crystal can realize a neuromorphic computer by encoding information in a low-intensity structured beam and exciting the crystal nonlinearity through a bias beam with rectangular spatial intensity distribution. Except for the presence of the structured beam, this model describes also the hydrodynamics of the dam break problem \citep{el2016dam}, tracing the path for further investigations in the field. \newline

When Eqs. (\ref{eq:RC}) are modified in a way that $\mathbf{x}$ represents spike neurons, we obtain an LSM \citep{maass2011liquid}, a NMC architecture proven to be a universal approximator. 
A hardware implementation of the LSM was realized at the beginning of this millennium \citep{fernando2003pattern}. This system uses mechanical excitation of linear waves in a bucket filled with water and, by imaging the wave interference at the center of the bucket, it accomplishes complex tasks such as vowel recognition or logic operations. Even if the understanding of this machine could clearly benefit from an analysis of its reservoir in terms of NSE, no research has been done in this direction because of the complexity of the system. However, we believe that such an investigation would be an important milestone to achieve in order to develop better models of the NMC's overall performance.

\section{Outlook to the NSE regularity problem and the soft matter research}
\label{outlook}
An answer to the NSE regularity problem will probably be the result of influx of ideas from other disciplines such as, theory of computation, theory of neural networks, and computational complexity. Fluid mechanics has been a topic of interest for engineers, physicists, mathematicians, biologists, and computer. The turn of 20th century saw advancements in the mathematical side when questions of existence and uniqueness of PDEs were asked, particularly by Ladyzhenskaya, Leray, and Oseen. There was a shift in thinking, in the mathematical physics community, from proving/finding a classical solution of a problem towards accepting the notion of a `weak solution' of the PDE. On the other hand, the experimental side of fluid mechanics, over the last few decades, has gained immense traction because of rapid developments in high resolution microscopy, high-speed imaging, experimentation in adverse conditions like space station, alongside the influences of molecular chemistry, polymer physics, quantum physics, and many other disciplines. The vast amounts of data available from experimental observations has invited an array of diverse neural network architectures to model the phenomenon. However, there still remains a lot to be done on the fundamental side of neural networks: there is a need for the theoretical-computational theory of fluid mechanics that can support ongoing remarkable feats of neural networks in modeling the experimental observations. In fact, a recent paper, Wang, Lai, Gomez-Serrano, Buckmaster \citep{wang2022self} showed that a physics-informed neural networks (PINN) can be employed to find smooth self-similar solutions for the Boussinesq equations. The authors argue that their numerical approach, introduced in their paper can guide a computer-assisted proof of finite-time singularity formation from finite-energy initial data. In fact, very recently, 2D Boussinesq and 3D Euler equations were proven to blow-up in self-similar fashion (\cite{chen2022stable}). The technique involved decomposing the solution into semi-analytical construction that captures the long-range behaviour and the numerical computation that yields the behaviour of the next time step, such that that combination of the two is useful to prove self-similar blowup of the solutions. This result present first of its own kind of computer-assisted proof of a blowup scenario. It is the hope that in the future, fluid equations without any boundary or support might well be proved to exhibit blowup for smooth initial data, however, that probably would involve highly computer-assisted programming of fluid parcels. \newline

Given that the Navier-Stokes regularity problem lies at the heart of current scientific advances of the 21st century, it is highly likely that the theoretical and numerical advances in physics-informed machine learning techniques and its connections to the computational complexity of `computable' fluid parcels are going to metamorphosise the physics of fluids, soft matter, and active matter in coming few decades and potentially lead to emergence of a new discipline along the lines of ``soft matter complexity, computing, and learning''. 

\section{Acknowledgments}
We thank Prof. Eva Miranda for additional references and points to add to this article; David Kaloper Meršinjak for long hours of discussion on various related problems that led to establishing the narrative of this article; Dr Shamit Shrivastava for the deeply insightful scientific conversations and for sharing with us his enthusiasm in understanding the complexity of the world. 

\subsection*{Dedication}
This article is dedicated to Sir George Gabriel Stokes, a physicist who made seminar contributions in the field of fluid mechanics. He served as the Lucasian Professor of Mathematics from 1849 until 1903, and served as the Master of Pembroke College from 1902-1903, at the University of Cambridge. The article written above is the result of never-ending source of inspiration that S.S. received by glancing fondly over the magnanimous portrait of Stokes in the dining hall of Pembroke.

% Bibliography
\bibliographystyle{unsrtnat}
\bibliography{reference}

\end{document}